\documentclass[12pt,]{article}

\usepackage{OldStandard}

\usepackage{amsmath,amssymb,amsfonts,amsthm}
\usepackage{hyperref}

\usepackage[sorting=ynt,backend=biber,style=alphabetic]{biblatex}

\newtheorem{theorem}{\scshape Theorem}

\title{\bfseries Exponent Lifting Property of Integer Sequences}
\author{\scshape Masum Billal}

\begin{document}
	\maketitle
		\begin{abstract}
			The purpose of this paper is to investigate integer sequences with exponent lifting property, a property common in Fibonacci or Lucas sequences.
		\end{abstract}
	An integer sequence $(a_{n})$ is a \textit{strong divisibility sequence} if $\gcd(a_{m},a_{n})=a_{\gcd(m,n)}$. Following \textcite{lehmer_1930,ward_1936}, the smallest index $\rho$ such that $p\mid a_{\rho}$ for a prime $p$ is the \textit{rank of apparition of $p$} in $(a_{n})$. Similarly, $\rho(p,r)$ is the rank of apparition of $p^{r}$ in $(a_{n})$. We say that $(a_{n})$ has the \textit{exponent lifting} property if for every prime divisor $p$ of $a_{n}$, we have
		\begin{align*}
			\nu_{p}(a_{nk})
				& = \nu_{p}(a_{n})+\nu_{p}(k)
		\end{align*}
	where $\nu_{p}(n)=k$ implies that $p^{k}\mid n,p^{k}\nmid n$. We can also denote it by $p^{k}\|n$. \textcite[Theorem $1$]{ward_1936} proves that $(a_{n})$ is a strong divisibility sequence if and only if for a prime $p$ and a positive integer $a$, $p^{a}\mid a_{k}$ if and only if $\rho(p,a)\mid k$. Note that the rank of apparition of $p$ in $(b_{n})$ is the same as the rank of apparition of $p$ in $(a_{n})$. Denote the product of the first $n$ terms of $(a_{n})$ by $n!_{a}$ (see \textcite[$\S3.3$]{billal_riasat_2021}).
		\begin{theorem}\label{thm:ltosd}
			If $(a_{n})$ has the exponent lifting property, then $(a_{n})$ is a strong divisibility sequence.
		\end{theorem}

		\begin{proof}
			For positive integers $m$ and $n$, let $g=\gcd(m, n),m=gu,n=gv$ where $\gcd(u,v)=1$ and $h=\gcd(a_{m},a_{n})$. We will show that $h=a_{g}$. First, consider that $p$ is a prime divisor of $g$. If $p^{e}\|a_{g}$,
				\begin{align*}
					\nu_{p}(h)
						& = \min\left(\nu_{p}(a_{gu}),\nu_{p}(a_{gv})\right)\\
						& = \nu_{p}(a_{g})+\min(\nu_{p}(u),\nu_{p}(v))
				\end{align*}
			Since $\gcd(u,v)=1$, $p$ cannot divide both $u$ and $v$. Therefore, either $\nu_{p}(u)$ or $\nu_{p}(v)$ is $0$ and $\min(\nu_{p}(u),\nu_{p}(v))=0$. This gives us $\nu_{p}(h)=\nu_{p}(a_{g})$ for all prime divisor $p$ of $g$. Next, assume that $p$ is a prime divisor of $h$ and $p^{e}\|h$. Then $p^{e}\mid a_{m}$ and $p^{e}\mid a_{n}$. More specifically, $p^{e}\|a_{gu}$ or $p^{e}\|a_{gv}$ must hold. Again, by definition $\nu_{p}(a_{gu})=\nu_{p}(a_{g})+\nu_{p}(u)$ and $\nu_{p}(a_{gv})=\nu_{p}(a_{g})+\nu_{p}(v)$. Since both $p\mid u$ and $p\mid v$ cannot hold, so $p^{e}\|a_{gu}$ or $p^{e}\|a_{gv}$ must hold. Then $p^{e}\|a_{g}$ holds for all $p^{e}\|h$. Thus, we must have $h=a_{g}$.
		\end{proof}
	If $(a_{n})$ is a strong divisibility sequence, then there is a sequence $(b_{n})$ such that
		\begin{align*}
			a_{n}
				& = \prod_{d\mid n}b_{d}
		\end{align*}
	where $\gcd(b_{m},b_{n})=1$ whenever $m\nmid n$ and $n\nmid m$, see \textcite[Chapter $3$]{billal_riasat_2021}. Such $(b_{n})$ is unique for $(a_{n})$ and is called the \textit{lcm sequence} of $(a_{n})$, also see \textcite{nowicki_2015}.
		\begin{theorem}
			Let $(a_{n})$ be a strong divisibility sequence, $(b_{n})$ be the lcm sequence of $(a_{n})$ and $\rho$ be the rank of apparition of prime $p$ in $(a_{n})$. Then $(a_{n})$ has the exponent lifting property if and only if for any positive integers $n$ and $m>1$ such that $p\nmid m$, $p\|b_{\rho p^n}$ but $p\nmid b_{\rho p^nm}$.
		\end{theorem}

		\begin{proof}
			First, we will prove the if part. Since $(a_{n})$ is a strong divisibility sequence, $p\mid a_{k}$ if and only if $\rho\mid k$. By assumption, $(a_{n})$ has exponent lifting property. If $p^\alpha\|a_\rho$, then $p^{\alpha+1}\|a_{\rho p}$.
				\begin{align*}
					a_{\rho p}
						& = \prod_{d\mid \rho p}b_d\\
					\nu_{p}(a_{\rho p})
						& = \nu_{p}\left(\prod_{d\mid\rho p}b_d\right)
				\end{align*}
			If $d<\rho$, then $p\nmid a_d$ so $p\nmid b_d$. Thus,
				\begin{align*}
					\nu_{p}(a_{\rho p})
						& = \nu_{p}\left(\prod_{d\mid p}b_{\rho d}\right)\\
						& = \nu_{p}(b_\rho)+\nu_{p}(b_{\rho p})\\
					\alpha+1
						& = \alpha+\nu_{p}(b_{\rho p})
				\end{align*}
			So, $\nu_{p}(a_{\rho p})=1$ and $p\mid b_{\rho p}$. By induction, we can see that $p$ not only divides $b_{\rho p^i}$ for $i\in\mathbf{N}$, more precisely, $p\|b_{\rho p^i}$. Next, assume that $p^{\alpha+u}\|a_{n}$ for some positive integer $n=\rho p^{u}m$ where $p\nmid m$. From the exponent lifting property and the argument above,
				\begin{align*}
					\nu_{p}(a_{n})
						& = \nu_{p}(a_{\rho p^{u}m})\\
						& = \nu_{p}(a_{\rho})+\nu_{p}\left(\prod_{d\mid p^um}b_{\rho d}\right)\\
						& = \alpha+\nu_{p}\left(\prod_{d\mid p^u}b_{\rho d}\right)+\nu_{p}\left(\prod_{\substack{d\mid p^u\\e\mid m\\e>1}}b_{\rho de}\right)\\
					\alpha+u
						& = \alpha+\sum_{i=1}^u\nu_{p}(b_{\rho p^i})+\nu_{p}\left(\prod_{i=1}^u\prod_{\substack{e\mid m\\e>1}}b_{\rho p^ie}\right)\\
						& = \alpha+u+\nu_{p}\left(\prod_{i=1}^u\prod_{\substack{e\mid m\\e>1}}b_{\rho p^ie}\right)\\
						& = \alpha+u+\sum_{i=1}^u\sum_{\substack{e\mid m\\e>1}}\nu_{p}(b_{\rho p^ie})
				\end{align*}
			From this, we have that $\nu_{p}(b_{\rho p^ie})=0$ for $1\leq i\leq u$ and $e\mid m$ if $e>1$. In other words, $p\mid b_k$ if and only if $k=\rho p^u$ for some non-negative integer $u$.

			For the only if part, we have that $(a_{n})$ is a strong divisibility sequence such that $p\| b_{\rho p^u}$ but $p\nmid b_{\rho p^um}$ for $m>1$. Let $n$ be a positive integer such that $n=\rho p^um$ and $p^\alpha\|a_\rho$.
				\begin{align*}
					\nu_{p}(a_{n})
						& = \nu_{p}(a_{\rho p^um})\\
						& = \nu_{p}\left(\prod_{d\mid \rho p^um}b_d\right)\\
						& = \nu_{p}(a_\rho)+\nu_{p}\left(\prod_{d\mid p^um}b_{\rho d}\right)\\
						& = \nu_{p}(a_{\rho})+\sum_{d\mid p^u}\nu_{p}(b_{\rho d})+\sum_{d\mid p^u}\sum_{\substack{e\mid m\\e>1}}\nu_{p}(b_{\rho de})\\
						& = \alpha+\sum_{i=1}^u\nu_{p}(b_{\rho p^i})+0\\
						& = \alpha+\sum_{i=1}^u1\\
						& = \alpha+u
				\end{align*}
			This proves the theorem.
		\end{proof}
	A corollary of this theorem is that all sequences with exponent lifting property are strong divisibility sequences but the converse may not hold. Another corollary is the following.
		\begin{theorem}
			Let $(a_{n})$ be a sequence with the exponent lifting property and $(b_{n})$ be the lcm sequence of $(a_{n})$. If $m$ and $n$ are distinct positive integers, then $\gcd(b_{m},b_{n})=1$ if and only if $m/n$ is a prime power.
		\end{theorem}
	To be more precise, if $p\mid b_{m}$ and $p\mid b_{n}$, then $m/n$ is a power of $p$. Next, we have an analogous of Legendre's theorem.
		\begin{theorem}\label{thm:leg}
			Let $(a_{n})$ be a sequence with the exponent lifting property. Then for a prime $p$,
				\begin{align*}
					\nu_{p}(n!_{a})
						& = \sum_{i\geq 1}\left\lfloor{\dfrac{n}{\rho(p,i)}}\right\rfloor
				\end{align*}
		\end{theorem}
	We see that the same argument as in the proof of Legendre's theorem applies here if we have that $\gcd(a_{m},a_{n})=a_{\gcd(m,n)}$. This condition is satisfied here by \autoref{thm:ltosd}. As a corollary of \autoref{thm:leg}, we have the next result.
		\begin{theorem}
			The binomial coefficients of an integer sequence with the exponent lifting property are integers.
		\end{theorem}

		\begin{proof}
			For any prime $p$,
				\begin{align*}
					\nu_{p}\left(\binom{n}{k}_{a}\right)
						& = \nu_{p}(n!_{a})-\nu_{p}(n!_{k})-\nu_{p}(n!_{n-k})\\
						& \sum_{i\geq 1}\left\lfloor{\dfrac{n}{\rho(p,i)}}\right\rfloor-\left\lfloor{\dfrac{k}{\rho(p,i)}}\right\rfloor-\left\lfloor{\dfrac{n-k}{\rho(p,i)}}\right\rfloor
				\end{align*}
			Since $\lfloor{x+y}\rfloor\geq\lfloor{x}\rfloor+\lfloor{y}\rfloor$, $\nu_{p}\left(\binom{n}{k}_{a}\right)\geq 0$ for all $p$. Thus, $\binom{n}{k}_{a}$ is an integer.
		\end{proof}
	\printbibliography
\end{document}